\documentclass[a4paper,12pt]{amsart}
\usepackage{amssymb}
\usepackage[]{latexsym}

\usepackage{amsthm}

\usepackage{lscape}

\usepackage[latin1]{inputenc}
\usepackage{color}
\usepackage[all,dvips]{xy}

\usepackage[dvips]{graphicx}

\newtheorem{thm}{Theorem}
\newtheorem{cor}[thm]{Corollary}

\newtheorem{prop}[thm]{Proposition}

\theoremstyle{remark}
\newtheorem{rem}[thm]{Remark}

\theoremstyle{definition}

\newenvironment{pf}{\par\noindent{\bf Proof.}\enspace\ignorespaces}{\qed\par\par}

\def\qed{\hfill $\Box$}

\usepackage[OT2,T1]{fontenc}
\DeclareSymbolFont{cyrletters}{OT2}{wncyr}{m}{n}
\DeclareMathSymbol{\Sha}{\mathalpha}{cyrletters}{"58}

\newcommand{\pmodd}[1]{\,({\rm mod\,}{#1}) }

\newcommand{\Q}{{\mathbb{Q}}}
\newcommand{\QD}{{\mathbb{Q}(\sqrt{D})}}
\newcommand{\Z}{{\mathbb{Z}}}
\newcommand{\ZZ}[1]{\mathbb{Z}/#1\mathbb{Z} }

\begin{document}

\title[Three cubes in arithmetic progression over quadratic fields]{Three cubes in arithmetic progression over quadratic fields}


\author[Enrique Gonz\'alez-Jim\'enez]{Enrique Gonz\'alez-Jim\'enez}
\address{Universidad Aut{\'o}noma de Madrid, Departamento de Matem{\'a}ticas and Instituto de Ciencias Matem{\'a}ticas (CSIC-UAM-UC3M-UCM), Madrid, Spain}
\email{enrique.gonzalez.jimenez@uam.es}

\thanks{The author was partially supported by the grant MTM2006--10548 and CCG08--UAM/ESP--3906.}

\date{\today}
\maketitle

\begin{abstract}
We study the problem of the existence of arithmetic progressions of three cubes over quadratic number fields $\QD$, where $D$ is a squarefree integer. For this purpose, we give a characterization in terms of $\QD$-rational points on the elliptic curve $E:y^2=x^3-27$. We compute the torsion subgroup of the Mordell-Weil group of this elliptic curve over $\QD$ and we give partial answers to the finiteness of the free part of $E(\QD)$. This last task will be translated to compute if the rank of the quadratic $D$-twist of the modular curve $X_0(36)$ is zero or not. 
\end{abstract}

\section{Introduction}

Nowadays, the study of arithmetic progressions consisting of perfect nth-powers is of considerable interest in number theory. Thanks to the development of new techniques to solve diophantine equations, several problems related to arithmetic progressions are being solved. For example, Darmon and Merel \cite{darmon-merel} proved that there are not non-trivial arithmetic progressions of three nth-powers. This article treated about the oldest and simplest problem in this direction, that is, the case of three cubes in arithmetic progression. According to Dickson's {\it History of the Theory of Numbers} \cite[ Vol. II, pp. 572--573]{dickson}, Legendre \cite{legendre} established that there are not non-trivial arithmetic progressions of three cubes over $\Q$. We are going to study when there exists a non-trivial arithmetic progression of three cubes over a quadratic number field. Note that over number fields, it has been obtained important progresses last years. In this direction, Xarles \cite{xarles} has proved that for any positive integers $n$ and $d$, the length of any arithmetic progression of nth-powers over a number field of degree $d$ is bounded by a constant depending only on $n$ and $d$. In particular, for the case of squares over quadratic fields, Xarles \cite{xarles} has proved that the length of any arithmetic progression of squares over any quadratic field is less than six. The case of length four and five has been treated in \cite{gonzalez-steuding} and \cite{gonzalez-xarles} respectively. Then the study of arithmetic progressions of squares over quadratic field is mainly done. Next step will be to study the case of arithmetic progression of cubes over quadratic fields. As a first step in this project, the goal of this paper is to study when there exists a three-term arithmetic progression consisting in cubes over $\QD$, where $D$ is a squarefree integer. For this purpose, first we will parametrize the arithmetic progressions of three cubes by the rational points of the elliptic curve $E:y^2=x^3-27$. Therefore to find three cubes in arithmetic progression over $\QD$ we should compute the Mordell-Weil group $E(\QD)$. Finally we will reduce our problem to the determination of the rank of the quadratic twists of the modular curve $X_0(36)$. We will use the work of Barthel \cite{barthel} and Frey \cite{frey} to obtain some partial answer to this question.

\section{Parametrization}
Let be $x_0^3,x_1^3,x_2^3$ three cubes in a field $k$. Assume that they form an arithmetic progression. Therefore, they satisfy $x_1^3-x_0^3=x_2^3-x_1^3$. That is,  the point $[x_0,x_1,x_2]\in\mathbb{P}^2(k)$ belongs to the projective curve $C:X_0^3-2X_1^3+X_2^3=0$. It is easy to check that if $\mbox{char($k$)}\neq 2,3$ then $C$ is an irreducible smooth projective curve of genus $1$ with two trivial points: $[1,1,1],[-1,0,1]\in C(k)$. Note that this two points correspond to the trivial arithmetic progressions: the constant progression $1,1,1$ and $-1,0,1$. Since the genus of $C$ is $1$ and $C$ has at least one rational point, $C$ is an elliptic curve defined over $k$. Let us compute a Weierstrass model for $C$. We have that  $[-1,0,1]$ is an inflection point of $C$. Let move the point $[-1,0,1]$ to $[0,1,0]$ and its tangent line to the line $w=0$. The tangent line at $[-1,0,1]$ is $X_0+X_2=0$, then the linear change of variables that sends $[X_0,X_1,X_2]$ to $[u,v,w]=[X_0,X_1,X_0+X_2]$ gives us the equation $-2v^3+3u^2w-3uw^2+w^3=0$. Now assuming that $\mbox{char($k$)}\neq 2,3$, we can make a change of variables to obtain an isomorphism to the elliptic curve $E:zy^2=x^3-27z^3$. This isomorphism is as follows:
$$
\varphi:C\longrightarrow E\,,\quad \varphi([x_0,x_1,x_2])=[6x_1,9(x_0-x_2),x_0+x_2]
$$
and its inverse is given by:
$$
\varphi^{-1}:E\longrightarrow C\,,\quad \varphi^{-1}([x,y,z])=\left[\frac{9z+y}{18},\frac{x}{6},\frac{9z-y}{18}\right].
$$

Therefore, we have proved the following proposition:

\begin{prop}\label{parametrization}
Let $k$ be a field of $\mbox{char($k$)}\neq 2,3$, then arithmetic progressions of three cubes in $k$ are parametrized by $k$-rational points of the elliptic curve $
E:zy^2=x^3-27z^3$. This parametrization is as follows:\\
\noindent $\bullet$ Let $x_0,x_1,x_2\in k$ such that $x_0^3,x_1^3,x_2^3$ form an arithmetic progression. Then $P=[6x_1,9(x_0-x_2),x_0+x_2]\in E(k)$.\\
\noindent $\bullet$ Let $P=[x,y,z]\in E(k)$. Define $x_0=9z+y,x_1=3x,x_2=9z-y$. Then $x_0^3,x_1^3,x_2^3$ form an arithmetic progression.
\end{prop} 

\begin{cor}\label{coro_para}
Let $k$ be a field of $\mbox{char($k$)}\neq 2,3$, then a necessary condition to the existence of a non-trivial arithmetic progressions of three cubes is the existence of a point $(x,y)\in E(k)$ such that $x\ne 3$. That is, $\ZZ{2}\subsetneq E(k)$.
\end{cor}

We will see that in general the condition $\ZZ{2}\subsetneq E(k)$ is not sufficient.

\

As a corollary we obtain:

\begin{cor}
There is not non-trivial arithmetic progressions of three rational cubes.
\end{cor}

\noindent This statement is due to Legendre \cite{legendre}. For the sake of completeness, we give the short using the above corollary.

\vspace{.3cm}

\begin{proof}  Using \verb+Sage+ \cite{sage} or \verb+Magma+ \cite{magma}, one can check that $E$ is the curve \verb+36A3+ in Cremona's tables \cite{cremona}, resp. \verb+36C+ in the Antwerp tables \cite{antwerp}. Checking these tables or using one of the above mentioned computer algebra systems, one can prove $E(\Q)\cong \ZZ{2}$. There are no $\Q$-rational affine points on $E$ apart $(3,0)$ which corresponds to the constant arithmetic progressions.
\end{proof}

\begin{rem}
Let $k$ a field, then if $\mbox{char($k$)}=2$ or $3$ looking for arithmetic progressions over $k$ is not interesting. If $\mbox{char($k$)}=2$ and $x_0^3,x_1^3,x_2^3$ is an arithmetic progression, then its length is $2$ instead of $3$, since $x_2^3=x_0^3$. Now, if $\mbox{char($k$)}=3$ then $C:X_0^3-2X_1^3+X_2^3=0$ is three copies of $X_0+X_1+X_2=0$, that is $C(k)\cong \mathbb{P}^1(k)$.  
\end{rem}

Our purpose in this paper is to obtain an answer to the question: are there non-constant arithmetic progressions of cubes over a quadratic number field? And in the affirmative case, give an explicit algorithm to construct them. Our main tool will be the characterization given at proposition \ref{parametrization}. 

Note that thanks to the above parametrization it is easy to check that for any $\alpha\in \Q$, we have an algebraic extension of $\Q$ adjoining to $\Q$ the squarefree part of $\alpha^3-27$ such that there exists a non-constant arithmetic progression of three cubes over that field. Nevertheless, this construction is not useful for our purpose, since we do not have control of the discriminant of this quadratic field. 

Therefore, for an squarefree integer $D$,  our goal is to compute the Mordell-Weil group of the elliptic curve $E:y^2=x^3-27$ over $\QD$. The torsion subgroup will be computed in section \ref{sec_torsion}. To compute the rank we will translated this problem to compute the rank of the quadratic $D$-twist of $E$ over $\Q$. This will be done in section \ref{sec_rank}.

\section{Torsion subgroup}\label{sec_torsion}

In this section we are going to give a complete characterization of the torsion subgroup of the elliptic curve $E:y^2=x^3-27$ over a quadratic number field $\QD$. We will denote by $E(\QD)_{\mbox{\tiny tors}}$ this subgroup. We can now prove the following result. 

\begin{prop}\label{prop_tors}
Let $D$ be a squarefree integer. Then the torsion subgroup of the elliptic curve $E:y^2=x^3-27$ over $\QD$ is
$$
E(\QD)_{\mbox{\tiny tors}}\cong\left\{
\begin{array}{lcl}
\ZZ{2} \oplus \ZZ{6}& & \mbox{if $D=-3$,}\\
\ZZ{2} & & \mbox{if $D\ne-3$.}
\end{array}
\right.
$$
\end{prop}

\begin{pf}
Kamienny \cite{kamienny} proved that the only primes possibly dividing the order of the torsion subgroup of an elliptic curve over a quadratic field are $2,3,5,7,11$ and $13$. Then it is enough to compute for which quadratic fields the elliptic curve $E:y^2=x^3-27$ has a torsion point of order $n\in\{2,3,4,5,7,11,13\}$. Note that we need to check $n=4$ since there are a point of order $2$ defined over $\Q$. To make this task we look for the irreducible factors of degree one or two of the $n$th-division polynomial of $E$ in $\Z[x]$. The set of these factors is $\{x,x-3,x^2 + 3 x + 9,x^2 - 6 x - 18\}$. Therefore the only possible values of $D$ such that  $E(\QD)_{\mbox{\tiny tors}}$ increases with respect $E(\Q)_{\mbox{\tiny tors}}$  are $D=3$ and $D=-3$. An straightforward computation shows that $E(\Q(\sqrt{3}))_{\mbox{\tiny tors}}\cong\ZZ{2}$ and $E(\Q(\sqrt{-3}))_{\mbox{\tiny tors}}\cong\ZZ{2}  \oplus \ZZ{6}$.
\end{pf}

\section{Rank}\label{sec_rank}
 The aim of this section is to compute the rank of the elliptic curve $E:y^2=x^3-27$ over a quadratic field. We are going to translate this problem to an easier one: to compute the rank of a quadratic twist of an elliptic curve over $\Q$.

\begin{prop}\label{prop_rank}
Let $D$ be a squarefree integer, $E:y^2=x^3-27$ and  $F^D:y^2=x^3+D^3$. Then 
$$
\mbox{rank}_{\,\Z}\, E(\QD)=\mbox{rank}_{\,\Z}\, F^D(\Q).
$$
\end{prop}

\begin{pf}
Let denote by $E^D$ the $D$-quadratic twist of $E$. That is, $E^D:y^2=x^3-27D^3$.  It is well known that for an arbitrary elliptic curve $E_0$ defined over $\Q$, we have  
\begin{equation}\label{eq_rank}
\mbox{rank}_{\,\Z}\, E_0(\QD)=\mbox{rank}_{\,\Z}\, E_0 (\Q) + \mbox{rank}_{\,\Z}\, E_0^D(\Q).
\end{equation}
Applying the above equality to $E$ we have $\mbox{rank}_{\,\Z}\, E(\QD)=\mbox{rank}_{\,\Z}\, E^D(\Q)$, since $E(\Q)$ is finite. 

Now, we have that $F^1:y^2=x^3+1$ is $\Q$-isogenous to $E^1=E$. This isogeny has the following equations
$$
 \psi:F^1\longrightarrow E^1\,,\qquad \psi(x,y)=\left(\frac{x^3+4}{x^2},\frac{x^3-8}{x^3}\, y\right).
$$
Therefore $F^D$ is $\Q$-isogenous to $E^D$, thus $\mbox{rank}_{\,\Z}\, F^D(\Q)=\mbox{rank}_{\,\Z}\, E^D(\Q)$. This finishes the proof.
\end{pf}

\

The study of the rank of the quadratic twists of an elliptic curve is an important area in the theory of elliptic curves. In particular, the quadratic twists of the elliptic curve $F^1$ have been deeply studied by Barthel \cite{barthel} and Frey \cite{frey}. Their results will be applied in the context of arithmetic progressions of three cubes in the next section.

\section{Arithmetic progressions of three cubes over quadratic fields}

\begin{thm}\label{teor}
Let $D$ be a squarefree integer. Then there is a non-trivial arithmetic progression of three cubes over $\QD$ if and only if the $D$-quadratic twist of $X_0(36)$ has positive rank.
\end{thm}

\begin{pf}
First we apply the characterization given at proposition \ref{parametrization} for the case $k=\QD$ obtaining that  arithmetic progressions of three cubes over $\QD$ are parametrized by $E(\QD)$. Corollary  \ref{coro_para} together with proposition \ref{prop_tors} tell us that the only possible $D$ such that there exists a non-trivial arithmetic progression of three cubes over $\QD$ coming from a torsion point of $E(\QD)$ is $D=-3$. Let $[x_0,x_1,x_2]\in \varphi^{-1} E(\Q(\sqrt{-3}))_{\mbox{\tiny tors}}$, then its corresponding arithmetic progression $x_0^3,x_1^3,x_2^3$ is equivalent to the arithmetic progression $-1,0,1$ or $1,1,1$. 

Now we are going to obtain non-torsion points on $E(\QD)$ coming from non-torsion points on $E^D(\Q)$. This will be done thanks to the following map
$$
 \phi:E^D\longrightarrow E\,,\qquad \phi(x,y)=\left(\frac{x}{D},\frac{y}{D^2}\sqrt{D}\right).
$$ 
Let $(x,y)\in E^D(\Q)$ then $\varphi^{-1}\circ\phi(x,y)=[9D^2-y\sqrt{D},3xD,9D^2+y\sqrt{D}]=[x_0,x_1,x_2]$ and denote by $\mathcal{S}$ the arithmetic progression $x_0^3,x_1^3,x_2^3$. First assume that $\mathcal{S}$ is equivalent to the arithmetic progression $-1,0,1$. Then $x=0$ and $y^2=-27D^3$, and since $y\in \Q$ we have $D=-3$ and $y=27$, that corresponds to the point $(0,3\sqrt{-3})\in E(\Q(\sqrt{-3}))[3]$. Now assume that  $\mathcal{S}$ is the constant arithmetic progression. Then we have $y=0$ since $(9D^2-y\sqrt{D})^3=(9D^2+y\sqrt{D})^3$. That is, $\mathcal{S}$ correspond to the point $(3,0)\in E(\Q)[2]$. Therefore we have proved that if $P\in E^D(\Q)$ is a non-torsion point, then $\varphi^{-1}\circ\phi(x,y)$ gives a non-trivial arithmetic progression of three cubes over $\QD$. To finish the proof just take into account that a Weierstrass model for $X_0(36)$ is $y^2=x^3+1$, therefore by proposition \ref{prop_rank} the proof is done.
\end{pf}

\begin{cor}\label{coro8}
Let $d$ be a squarefree positive integer coprime with $6$ and 
$$
A_d=\!\!\!\!\!\sum_{m,n,k\in S}(-1)^n
\,\,\mbox{where}\,\,
S=\left\{
 m,n,k\in\Z \,\Big|\,
\begin{array}{l}
m^2+n^2+k^2=d,\\
m\equiv 1 \pmodd{3},n\equiv 0 \pmodd{3}\\
m+n\equiv 1 \pmodd{2}
\end{array}
\!\!\! \right\}.
$$
{\it (a)} Assuming the Birch \& Swinnerton-Dyer Conjecture: if $A_d = 0$ then there is a non-trivial arithmetic progression of three cubes over $\Q(\sqrt{-d})$.\\ 

\noindent {\it (b)} If $A_d \ne 0$ then there is not a non-trivial arithmetic progression of three cubes over $\Q(\sqrt{-d})$.\\

\end{cor}

\begin{pf}
Barthel \cite{barthel} and Frey \cite{frey} independently found a modular form $\Phi\in S_{3/2}(144,1)$ such that its image by the Shimura correspondence is the modular form $f\in S_2(36,1)$ attached to the elliptic curve $F:y^2=x^3+1$. Note that $F=X_0(36)$. That is, if we denote by $\mathcal{S}h$ the Shimura correspondence \cite{shimura} that maps a weight $3/2$ modular form to a weight $2$ modular form then $\mathcal{S}h(\Phi)=f$. Now, the $q$-expansion of $\Phi$ is
$$
\Phi(q)=\sum_{n\geq 1}A_nq^n.
$$
Applying Waldspurger's results \cite{waldspurger} to the elliptic curve $F$, they proof that if $d$ is a squarefree positive integer coprime with $6$ then 
$$
L(F^{-d},1)=0 \quad \mbox{if and only if}\quad A_d=0. 
$$
Therefore, if $A_d\ne 0$ we have that $L(F^{-d},1)\ne 0$ and by Kolyvagin \cite{koly} the rank of $F^{-d}(\Q)$ is equal to zero. This proof (b), by Theorem \ref{teor}. Assuming the Birch \& Swinnerton-Dyer Conjecture  it follows that if $A_d=0$ then $F^{-d}(\Q)$ is infinite. Again Theorem \ref{teor} finishes the proof of (a).
\end{pf}

\begin{cor}
Let $D$ be a squarefree integer and $\varepsilon\in\{\pm 1\}$. \\
{\it (a)} There is a non-trivial arithmetic progression of three cubes over $\QD$ if:
\begin{itemize}
\item[(i)] $D=\varepsilon p$ where $p>3$ is a prime such that $p\equiv 3\pmodd{4}$.
\item[(ii)] Assuming the Birch \& Swinnerton-Dyer Conjecture:
\begin{itemize}
\item[$\bullet$] $D>0$ and $D$ even coprime with $3$.
\item[$\bullet$] $D<0$ and $D\equiv 1, 5\pmodd{12}$.
\end{itemize}
\end{itemize}

\noindent {\it (b)}  There is not a non-trivial arithmetic progression of three cubes over $\QD$ if: 
\begin{itemize}
\item[(i)] $D$ such that if a prime $p$ divides $D$ then  $p\equiv 5\pmodd{12}$ or $p=3$.
\item[(ii)] $D=-p$ where $p$ is a prime such that $p\equiv 1\pmodd{12}$ and $x^4+3=0$ has not solution over $\mathbb{F}_p$. 
\end{itemize}
\end{cor}

\begin{pf}
This corollary is basically a translation of the results of Barthel \cite{barthel} and Frey \cite{frey} on the study of the rank of the $D$-quadratic twist of the elliptic curve $y^2=x^3+1$ to our context using the Theorem \ref{teor}. Note that Barthel only treated the case of $D$ negative and he only used Waldspurger's results and Shimura's correspondence {\it a la} Tunnell \cite{tunnell} to obtain her results. Meanwhile, Frey treated also the positive case. He used several techniques like Heegner points, $2$-descent and the above method used by Barthel too. 

Frey \cite[Prop. 5]{frey} proved that if $p$ is a prime greater than $3$ such that $p\equiv 3\pmodd{4}$ then $\mbox{rank}_{\,\Z}\, F^{\varepsilon p}(\Q)=1$. This implies (a)(i). 

Barthel showed that the functional equation of $L(F^D,s)$ satisfies that $L(F^D,1)=0$ if $D>0$ and $D$ even coprime with $3$ or $D<0$ and $D\equiv 1, 5\pmodd{12}$. Therefore assuming the Birch \& Swinnerton-Dyer Conjecture we have that for the above values of $D$, the rank of $F^D(\Q)$ is positive. Therefore the proof of (a)(ii) is finished.

 By \cite[Prop. 3 \& Bemerkungen p.82]{frey} we have that if all the prime divisors of $D$ are $5$ modulo $12$ then the rank of $F^{D}(\Q)$ is zero. Now let $D$ be with the above condition. Then applying the equality (\ref{eq_rank}) to the elliptic curve $F^{-3}$ and $-D$ and taking into account that $F^{-3}$ is $\Q$-isogenous to $F^1$ we have that 
$$
 \mbox{rank}_{\,\Z}\, F^{3D}(\Q) = \mbox{rank}_{\,\Z}\, F^{-D}(\Q)=0,
$$
that proves (b)(i). 

Finally, if $p$ is a prime such that $p\equiv 1\pmodd{12}$ and $x^4+3=0$ has not solution over $\mathbb{F}_p$ then $A_p\ne 0$ (cf. \cite[Prop. 2]{barthel} or \cite[Kor. 2]{frey}) and then by Corollary \ref{coro8} the proof of (b)(ii) is finished.
\end{pf}

\

{\it Remark.}  Frey  \cite[Satz 4, p. 73]{frey} affirms that  if $p$ is a prime such that $p\equiv 1\pmodd{12}$ and it is not completely split over $\Q(\sqrt[4]{\varepsilon 3})$, then $\mbox{rank}_{\,\Z}\, F^{\varepsilon p}(\Q)= 0$. He proved the case $\varepsilon=-1$ at  Korollar 2. But the case $\varepsilon=1$ is not true. For example, for $p=37$ we have that $\mbox{rank}_{\,\Z}\, F^{37}(\Q)= 2$ and $37\mathcal{O}=\mathfrak{p}_1\mathfrak{p}_2$ is the ideal prime factorization, where $\mathcal{O}$ is the ring of integer of $\Q(\sqrt[4]{3})$.

\subsection{Computational results}

Using the functionality \verb+mwrank+ on \verb+Sage+, we may compute the rank of $E^D(\Q)$; if this rank is non-zero, we can also compute an explicit arithmetic progression of three cubes over $\QD$: Let $P=(x,y)\in E^D(\Q)$ of infinite order, then $(9 D^2+y\sqrt{D})^3,(3 x D)^3,(9 D^2-y\sqrt{D})^3$ is an arithmetic progression over $\QD$. The following table list explicit examples of such progressions for the range $|D|\le 30$. At the first column  indicates the value of $D$ and the second gives a point $P=(x,y)\in E^D(\Q)$ of infinite order.

\medskip

\begin{center}
\begin{tabular}{|c||c|}
\hline
$D$ & $P=(x,y)\in E^D(\Q)$ with $\mbox{ord}\,P=\infty$\\
\hline
$-30$ & $ (-54, 756)$\\
\hline
$-26$ & $ (-26, 676)$\\
\hline
$-23$ & $ (987505/24336, -2386987127/3796416)$\\
\hline
$-21$ & $ (189, 2646)$\\
\hline
$-19$ & $ (-38, 361)$\\
\hline
$-11$ & $ (-6, 189)$\\
\hline
$-7$ & $ (7, 98)$\\
\hline
$-6$ & $ (9, 81)$\\
\hline
$2$ & $ (10, 28)$\\
\hline
$7$ & $ (1785/4, 75411/8)$\\
\hline
$10$ & $ (946/9, 28756/27)$\\
\hline
$11$ & $ (178849/400, -75621007/8000)$\\
\hline
$14$ & $ (217, 3185)$\\
\hline
$19$ & $ (1173649/2025, 1270868732/91125)$\\
\hline
$21$ & $ (126, -1323)$\\
\hline
$22$ & $ (22825/36, -3446443/216)$\\
\hline
$23$ & $ (4655599441/56851600, -201357032252761/428661064000)$\\
\hline
$26$ & $ (28249/100, 4697693/1000)$\\
\hline
\end{tabular}
\end{center}

\vspace{.3cm}

{\it Example:} Let $P=(10, 28)$ be a generator of the free part of the Mordell-Weil group $E^2(\Q)$.  The morphism  $\phi:E^2\rightarrow E$ applied to the point $P$ gives 
$$\phi(10,28)=\left(5,7\sqrt{2}\right)\in E(\Q(\sqrt{2})).$$
Now, the isomorphism $\varphi^{-1}:E\rightarrow C$ gives 
$$\varphi^{-1}([5,7\sqrt{2},1])=\left[\frac{9+7\sqrt{2}}{18},\frac{5}{6},\frac{9-7\sqrt{2}}{18}\right]\in C(\Q(\sqrt{2})),$$ that corresponds to the arithmetic progression $(9+7\sqrt{2})^3, (15)^3, (9-7\sqrt{2})^3$ over $\Q(\sqrt{2})$.

\end{document}